\documentstyle[12pt]{amsart}  % amslatex 1.1

\newtheorem{thm}{Theorem}[subsection]

\newtheorem{cor}[thm]{Corollary}
\newtheorem{prop}[thm]{Proposition}
\newtheorem{lemma}[thm]{Lemma}

\theoremstyle{remark}
\newtheorem{remark}[thm]{Remark}

\theoremstyle{definition}
\newtheorem{defn}[thm]{Definition}

\numberwithin{equation}{section}

\newcommand{\shTor}{{\underline{\operatorname{Tor}}}}
\newcommand{\Tor}{{\operatorname{Tor}}}
\newcommand{\Spec}{\operatorname{Spec}}

\newcommand{\supp}{\operatorname{Supp}}
\newcommand{\Supp}{\operatorname{Supp}}

\newcommand{\codim}{\operatorname{codim}}
\newcommand{\Ch}{\operatorname{Ch}}
\newcommand{\Ext}{\operatorname{Ext}}
\newcommand{\depth}{\operatorname{depth}}

\newcommand{\DR}{\operatorname{DR}}
\newcommand{\aDR}{\operatorname{DR}^{-1}}

\newcommand{\Fl}{F^\bullet}

\newcommand{\Gr}{\operatorname{Gr}}

\newcommand{\bbE}{{\Bbb E}}
\newcommand{\bbC}{{\Bbb C}}

\newcommand{\bbD}{{\Bbb D}}
\newcommand{\bbQ}{{\Bbb Q}}
\newcommand{\bbL}{{\Bbb L}}

\newcommand{\cF}{{\cal F}}
\newcommand{\cO}{{\cal O}}
\newcommand{\cM}{{\cal M}}

\newcommand{\cD}{{\cal D}}

\newcommand{\cU}{{\cal U}}

\newcommand{\m}{{\frak m}}

\newcommand{\dl}{d}

\newcommand{\dbar}{\overline{\partial}}

\newcommand{\xcirc}{X^\circ}

\newcommand{\isomo}{\overset{\sim}{=}}

\newcommand{\ltcx}{{\cal A}^{\bullet}_{(2)}({\Bbb E})}

\newcommand{\R}{\operatorname{\bold R}}

\newcommand{\Hom}{\operatorname{Hom}}

\newcommand{\lt}{$L^2$}

\newcommand{\shHom}{\underline{\operatorname{Hom}}}

\begin{document}

\title{Filtered perverse complexes}
\author{Paul Bressler}
\address{Department of Mathematics\\Pennsylvania State University\\
University Park, PA, USA}
\email{bressler@@math.psu.edu}
\author{Morihiko Saito}
\address{Research Institute for Mathematical Sciences\\
Kyoto University\\Kyoto, Japan}
\author{Boris Youssin}
\address{Department of Mathematics and  Computer Science\\
University of the Negev\\Be'er Sheva', Israel}
\email{youssin@@cs.bgu.ac.il}
\thanks{The last author was partially
supported  by NSF grants at  MIT and Indiana University, Bloomington}
\maketitle

\begin{abstract}
We introduce the notion of {\em filtered perversity} of a filtered differential
complex on a complex analytic manifold $X$, without any assumptions of
coherence, with the purpose of studying the connection between
the pure Hodge modules and the \lt-complexes.

We show that if a filtered differential complex $(\cM^\bullet,F_\bullet)$
is filtered
perverse then $\aDR(\cM^\bullet,F_\bullet)$ is isomorphic to a filtered
$\cD$-module; a coherence assumption on the cohomology of
$(\cM^\bullet,F_\bullet)$ implies that, in addition, this $\cD$-module
is holonomic.

We show the converse: the de Rham complex of a holonomic Cohen-Macaulay
filtered $\cD$-module is filtered perverse.
\end{abstract}

\section{Introduction}

\subsection{Cheeger---Goresky---MacPherson conjectures}

J.~Cheeger, M.~Goresky and R.~MacPherson \cite{CGM} conjectured some fifteen
years ago that
the intersection cohomology of a singular complex projective algebraic variety
is naturally isomorphic to its $L^2$ cohomology and the K\"ahler
package holds for them.
Their motivation was as follows.

The intersection cohomology was discovered by M.~Goresky and R.~MacPherson
\cite{GM1}, \cite{GM2} as an invariant of stratified spaces which for
complex algebraic varieties might serve as a replacement of the usual
cohomology: it had some properties that
the usual cohomology of smooth projective varieties possessed but
the usual cohomology of singular projective varieties did not.
One of such properties was Poincar\'e duality which is a part of the ``K\"ahler
package'' of properties that hold in the smooth case.

At the same time, J.~Cheeger discovered that the $L^2$ cohomology
groups of varieties with conical singularities have properties similar
to those of intersection cohomology, and he proved in this case
the Hodge---de Rham isomorphism between the $L^2$ cohomology that he defined
and studied, and the intersection cohomology \cite{C}.

The hope that underlied these conjectures was that it would be possible to
use the $L^2$ K\"ahler methods to prove the K\"ahler package for intersection
cohomology similarly to the way  the K\"ahler package was proved for the
usual cohomology of complex projective manifolds.

The most important part of the K\"ahler package is the $(p,q)$-decomposition
in the cohomology groups (the ``Hodge structure'').

The definition of $L^2$ cohomology involves a metric
(Riemannian or K\"ahler) defined almost everywhere on the variety (e.g.\ on
its nonsingular part).  The most important metric comes from a projective
embedding of the variety and is induced by the Fubini---Studi metric on
the projective space.  (The $L^2$ cohomology is independent of the choice of
the imbedding.)

The isomorphism with intersection cohomology is known in case of surfaces
\cite{HP}, \cite{Nag1} and in case of isolated singularities of any dimension
both for Fubini---Studi metric \cite{O2}, \cite{O2a} and for a different,
complete metric, introduced by L.~Saper,
which is defined on the nonsingular part of the variety
and blows up near the singularities \cite{Sap}.
The $(p,q)$-decomposition is known for the case of Fubini---Studi metrics
only in cases of dimension two \cite{Nag2} (except for the middle degree
cohomology groups) while a classical
result of Andreotti---Vesentini implies the $(p,q)$-decomposition
for {\em any\/} complete metric.
The general case is still open, despite the announcement of T.~Ohsawa
\cite{O3}.

In the meantime the second author \cite{S1}, \cite{S2} developed a theory of
polarizable Hodge modules which implied the K\"ahler package
for the intersection cohomology.
His main tool was the theory of $\cD$-modules and his methods were essentially
algebraic, reducing the intersection cohomology to the intersection cohomology
of a curve with coefficients
in a polarised variation of Hodge structure \cite{Z1}.

\subsection{The comparison between the Hodge structures}

Assuming that the Cheeger---Goresky---MacPherson
conjectures are true, one is faced with the question of comparison between the
two Hodge structures on the intersection cohomology: one induced by the
isomorphism with {\lt} cohomology, the other coming from the theory of
polarised Hodge modules.
In fact, different metrics give different $L^2$ cohomology theories and hence,
pose different comparison problems.

In case of isolated singularities, S.~Zucker \cite{Z} proved the coincidence
between the Hodge structures coming from polarized Hodge modules and from
$L^2$ cohomology with respect to the
Saper metric (or arithmetic quotient metrics similar to it).
Some partial results are also known in case of Fubini---Studi metric, see
\cite{Z} and \cite{Nag2}.

It is  interesting to note that the original purpose of the conjectures was
to construct the Hodge structure on the intersection cohomology.
The $L^2$ methods, however, turned out to be so difficult
that the Hodge structure was constructed by different,
algebraic methods and now we
are faced with the problem of comparison between the two Hodge structures.

\subsection{The local comparison problem}

A major component of a polarizable Hodge module is a regular
holonomic $\cD$-module $M$ with a good filtration $F_\bullet$.
Suppose $M$ corresponds to the intersection cohomology
complex of a complex projective subvariety $Z$;
the correspondence is given by taking the de Rham complex $\DR(M)$ of $M$
(so that $\DR(M)$ is isomorphic to the intersection cohomology complex of $Z$).
Then the filtration $F_\bullet$ induces a filtration on $\DR(M)$
which yields the Hodge structure on the intersection cohomology.

The complex $\DR(M,F_\bullet)$ is a filtered differential complex~\cite{S1}:
a complex of sheaves which are modules over the sheaf of analytic functions
and the differentials are differential operators.
This filtered differential complex completely determines $(M,F_\bullet)$ as
there is an inverse functor $\aDR$ \cite{S1}.

If the metric used in the construction of the {\lt} complex
(it is a K\"ahler metric on the nonsingular
part of $Z$) is bounded below with respect to Fubini---Studi metric,
the {\lt} complex is a filtered differential complex (see, e.g.,~\ref{subsec:l2}
below).

The {\em local comparison problem\/} is as follows: is it true that
the de Rham complex of the Hodge module is isomorphic to the {\lt} complex
in the derived category of filtered differential complexes?

The intersection cohomology can be taken with coefficients in a local system
defined on the non-singular part of $Z$ or a Zariski-open subset of it.
If this local system underlies a polarized variation of Hodge
structures then a corresponding polarized Hodge module can be constructed
and the intersection cohomology with coefficients in this local system has
a Hodge structure.
On the other hand, the $L^2$ cohomology can be taken with coefficients in
the same polarized variation, and we can ask the same local comparison
question in this situation.

\subsection{Weak filtered perversity}

A way to approach this problem is to try to identify the properties of
a filtered differential complex $(\cM^\bullet,F_\bullet)$ which
would imply that $\aDR(\cM^\bullet,F_\bullet)$ is isomorphic to the filtered
$\cD$-module which underlies a polarized Hodge module.

In general, $\aDR(\cM^\bullet,F_\bullet)$ is a
{\em complex\/} of filtered $\cD$-modules; in this paper we study
properties of $(\cM^\bullet,F_\bullet)$ which imply that this complex
is isomorphic to one filtered $\cD$-module in the filtered derived
category.

We call these properties {\em weak filtered perversity\/} (see
Definition~\ref{cond:wfp}); it means that, first, the complex
$(\cM^\bullet,F_\bullet)$ is locally trivial along the strata --- in a certain
filtered sense --- with respect to some analytic stratification, and second,
it satisfies certain local filtered cohomology vanishing which is similar
to the local cohomology vanishing of perverse sheaves.
No coherence assumption is being made on $(\cM^\bullet,F_\bullet)$.

In case $(\cM^\bullet,F_\bullet)$ is the \lt\ complex, the cohomology that must
vanish, turn out to be a version of the \lt-$\dbar$-cohomology,
see~\S\ref{rem:bound-cond} below for the discussion.

\subsection{The main results}

We show (see Theorem~\ref{thm:main})
that if $(\cM^\bullet,F_\bullet)$ is weakly filtered perverse then,
indeed, $\aDR(\cM^\bullet,F_\bullet)$ is isomorphic to a filtered $\cD$-module.

We show the converse (see Theorem~\ref{thm:converse}):
if $(M^\bullet,F_\bullet)$ is a coherent filtered $\cD$-module
which is Cohen-Macaulay (i.e., its dual in the filtered sense is also a
complex of filtered $\cD$-modules isomorphic to one filtered $\cD$-module)
then $\DR(M,F_\bullet)$ is weakly filtered perverse.

We show (see Proposition~\ref{cor:holon}) that a coherence assumption
together with filtered perversity of $(\cM^\bullet,F_\bullet)$ implies
that the filtered $\cD$-module $\aDR(\cM^\bullet,F_\bullet)$ is holonomic.

\subsection{Plan of the paper}

In Section~\ref{sec:fdm-dc} we review the necessary background material
from~\cite{S1}.

In Section~\ref{sec:fil.perv} we introduce the notion of weak filtered
perversity.

In Sections \ref{sec:fpc-to-dm} and~\ref{sec:coh.case} we prove the results
listed above.

In Section~\ref{sec:appl} we give a modest application: we strengthen the
results of \cite{KK2} and \cite{S1} and show (in the situation of \cite{KK2})
that filtered perversity of the $L^2$ complex implies
the local filtered isomorphism
(in the sense of derived category) between the $L^2$ complex
and the de Rham complex of
the $\cD$-module that underlies the corresponding pure Hodge module.

\subsection{Acknowledgements}
It is our pleasant duty to express our heartful thanks to all people who
helped us with their advice and helpful discussions: Daniel Barlet,
Alexander Beilinson, Joseph Bernstein, Jean-Luc Brylinski,
Michael Kapranov, Masaki Kashiwara, David Kazhdan,
Takeo Ohsawa, Claude Sabbah.

\section{Filtered $\cD$-modules and differential complexes}
\label{sec:fdm-dc}

In this section we
make a brief survey of the necessary parts of \cite{S1}.
%and refer the reader there for further details.

\subsection{General notation}
Throughout this paper $X$ will denote a complex manifold, $\cO_X$ the sheaf of
holomorphic functions, $\omega_X$ the canonical sheaf of $X$,
$\cD_X$ the sheaf of differential operators.
Unless specified otherwise, a $\cD_X$-module will always refer to a
sheaf of {\em right} modules over $\cD_X$.

For a complex of sheaves $\cF^\bullet$, we shall denote by
$H^j\cF^\bullet$ its {\em sheaf\/} cohomology, and by
$H^j_{\{x\}}\cF^\bullet$ its (hyper)cohomology with supports in a one-point
set $\{x\}$.

\subsection{Filtered $\cD$-modules}
Recall that $\cD_X$ is a filtered ring when equipped with the filtration
$F_\bullet\cD_X$, where $F_p\cD_X$ is the $\cO_X$-module of $\cD_X$ of
operators of order at most $p$.

A filtered $\cD_X$-module is a pair $(M,F_\bullet)$ consisting of a
$\cD_X$-module $M$ and a filtration $F_\bullet M$ of $M$ by $\cO_X$
submodules compatible with the action of $\cD_X$ and the filtration
on the latter.

We refer the reader to~\cite{S1}, \S2.1 for the precise definition of the
category of filtered $\cD_X$-modules and its derived category, in various
flavors; what is important for us now, is that the derived category
of filtered $\cD_X$-modules is
isomorphic to the derived category of the category whose objects are filtered
$\cO_X$-modules and whose morphisms are differential operators
that agree with the filtration in a certain way.

This equivalence is given by the two functors, $\DR_X$ and $\aDR_X$ (loc.~cit.)
which act as follows.

\subsection{The de Rham functor $\protect\DR_X$}
For a filtered $\cD_X$-module $(M,F_\bullet)$, the filtered differential
complex $\DR_X(M,F_\bullet)$ is the usual de Rham complex of $M$, given by
\[
\DR_XM=
M\otimes_{\cD_X}\left(\cD_X\otimes_{\cO_X}
\textstyle\bigwedge^{-\bullet}\Theta_X\right)
=M\otimes_{\cO_X}\textstyle\bigwedge^{-\bullet}\Theta_X
\]
where $\Theta_X$ is the tangent sheaf to $X$,
$\bigwedge^{-\bullet}\Theta_X$ its exterior algebra
with $p$-th exterior power placed in degree $-p$,
the differential is given by
\begin{equation*}
\begin{split}
{\dl}(P\otimes\xi_1\wedge\ldots\wedge\xi_p) & =\sum_{i=1}^{p}
(-1)^{i-1}P\xi_i\otimes\xi_1\wedge\ldots\wedge\widehat{\xi_i}
\wedge\ldots\wedge\xi_p \\
& +\sum_{1\leq i < j\leq p}(-1)^{i+j}P\otimes
\lbrack\xi_i,\xi_j\rbrack\wedge\xi_1\wedge\ldots\wedge\widehat{\xi_i}
\wedge\ldots\wedge\widehat{\xi_j}\wedge\ldots\wedge\xi_p
\end{split}
\end{equation*}
(it corresponds to the differential in
$\cD_X\otimes_{\cO_X}\bigwedge^{-\bullet}\Theta_X$
which makes it the standard
Koszul resolution of $\cO_X$ as a $\cD_X$-module),
and the filtration on $\DR_XM$ is given by
\[
F_p\left(M\otimes_{\cO_X}\textstyle\bigwedge^{-i}\Theta_X\right)=
F_{p+i}M\otimes_{\cO_X}\textstyle\bigwedge^{-i}\Theta_X\ .
\]
For a complex of filtered $\cD_X$-modules $(M^\bullet,F_\bullet)$, the
filtered differential complex $\DR_X(M^\bullet,F_\bullet)$ is the total
complex of $\DR_X(M^q,F_\bullet)$ for all $q$.

(Note that what we denote by $\DR_X$, was denoted by
$\widetilde\DR$ in~\cite{S1}.)

\subsection{The inverse de Rham functor $\protect\aDR_X$}
For a filtered differential complex $(\cM^\bullet,F_\bullet)$,
the complex of filtered $\cD_X$-modules $\aDR_X(\cM^\bullet,F_\bullet)$
is described as the complex of differential operators from $\cO_X$
into $(\cM^\bullet,F_\bullet)$ with the obvious differential and filtration.
The action of $\cD_X$, i.~e., differential operators $\cO_X\to\cO_X$,
is by composition.
The individual terms of this complex are simply $\cM^j\otimes_{\cO_X}\cD_X$.

The two functors $\DR_X$ and $\aDR_X$ are inverse to each other in
the derived categories.

\subsection{Duality}
\label{subsec:duality}
For a bounded complex of filtered  $\cD_X$-modules
$(M^\bullet,F_\bullet)$, its dual
$\bbD(M,F_\bullet)$ is is another complex of filtered  $\cD_X$-modules
defined (\cite{S1}, 2.4.3)
%as\[\shHom^\bullet_{\cD_X}(M^\bullet,\text{ filtered dualizing complex})\]
in such way that it agrees with various other duality functors, as follows.

There is a duality functor (also denoted by $\bbD$) on filtered differential
complexes which on an individual $\cO_X$-module $L$ is defined as
$\bbD(L)=\R\shHom^\bullet_{\cO_X}(L,\omega_X[n])$ where $n=\dim X$
(\cite{S1}, 2.4.11),
and such that in the appropriate derived categories the functors
$\bbD\circ\DR_X$ and $\DR_X\circ\bbD$ are isomorphic.

For a filtered differential complex $(\cM^\bullet,F_\bullet)$ we have
\begin{equation}
\label{eqn:D.vs.Gr}
\Gr^F_\bullet\bbD(\cM^\bullet,F_\bullet)\isomo
\R\shHom^\bullet_{\cO_X}(\Gr^F_\bullet\cM^\bullet,\omega_X[n])\ .
\end{equation}

In case $(M^\bullet,F_\bullet)$ is a complex of filtered coherent
$\cD_X$-modules (i.e., coherent $\cD_X$-modules
with good filtrations), the complex of $\cD_X$-modules which underlies
$\bbD(M^\bullet,F_\bullet)$, is the usual dual of $M^\bullet$.
In addition, under the same assumptions
$\bbD\bbD(M^\bullet,F_\bullet)\isomo(M^\bullet,F_\bullet)$.

\subsection{Restriction to a noncharacteristic submanifold}
Let $Y$ be a smooth submanifold of codimension $d$ in $X$; denote
the embedding $i:Y\hookrightarrow X$ and
$\omega_{Y/X}=\omega_Y\otimes_{\cO_X}\omega_X^{-1}$.

We say that
a bounded complex of filtered $\cD_X$-modules $(M^\bullet,F_\bullet)$
is {\em weakly noncharacteristic\/} with respect to $Y$ (or $Y$ with
respect to $(M^\bullet,F_\bullet)$) if it satisfies the property
\begin{equation}
\label{eqn:tor-nonchar}
\shTor^{\cO_X}_k(H^j(\Gr^F_p M^\bullet),\cO_Y) = 0
\text{ \ \ for all $k\ne0$, $j$ and $p$.}
\end{equation}

Under this assumption, the noncharacteristic restriction
$(M^\bullet,F_\bullet)_Y$ is defined as follows:
\[
(M^\bullet,F_\bullet)_Y=
(M^\bullet,F_\bullet)\otimes^\bbL_{\cD_X}(\cD_{X\leftarrow Y},F_\bullet)
\]
where
$\cD_{X\leftarrow Y}=\cD_X\otimes_{\cO_X}\omega_{Y/X}$
is the usual $(\cD_X,\cD_Y)$-bimodule with the filtration
$F_p\cD_{X\leftarrow Y}=
F_{p-d}\cD_X\otimes_{\cO_X}\omega_{Y/X}$;
the restriction $(M^\bullet,F_\bullet)_Y$ thus defined,
is a complex of right $\cD_Y$-modules.
As a complex of $\cO_Y$-modules, it can be described as
$(M^\bullet)_Y=M^\bullet\otimes^\bbL_{\cO_X}\omega_{Y/X}$ with
the filtration
$F_p(M^\bullet)_Y=F_{p-d}M^\bullet\otimes^\bbL_{\cO_X}\omega_{Y/X}$.
We have
\begin{equation}
\label{eqn:HGr(restr)}
H^j\Gr^F_p\left((M^\bullet,F_\bullet)_Y\right)\isomo
H^j(\Gr^F_{p-d} M^\bullet)\otimes_{\cO_X}\omega_{Y/X}\ .
\end{equation}

Suppose that $(M^\bullet,F_\bullet)$ has the property that the complex
$\Gr^F_\bullet M^\bullet$ has bounded $\Gr^F_\bullet\cD_X$-coherent
cohomology; in such case we say that $(M^\bullet,F_\bullet)$
is {\em noncharacteristic\/} with respect to $Y$ if, first,
\eqref{eqn:tor-nonchar} is satisfied, and second,
$\Gr^F_\bullet (M^\bullet)_Y$ also has bounded $\Gr^F_\bullet\cD_Y$-coherent
cohomology.

In the particular case when the complex $(M^\bullet,F_\bullet)$
is actually a filtered coherent $\cD_X$-module $(M,F_\bullet)$,
this definition is equivalent to the definition in \cite{S1}, 3.5.1 because
the condition of coherence of $(M^\bullet,F_\bullet)_Y$ is is equivalent
to the finiteness of the projection
$(Y\times_X T^\ast X)\cap\Ch(M)\to T^\ast Y$ where $\Ch(M)$ denotes the
characteristic variety of $M$.
In such case if $Y$ is noncharacteristic, we have
$i^\ast(M,F_\bullet)=(M^\bullet,F_\bullet)_Y[d]$ and
$i^!(M,F_\bullet)$ is isomorphic to $(M^\bullet,F_\bullet)_Y[-d]$ up to a
shift of filtration.

\begin{defn}
A filtered coherent $\cD_X$-module $(M,F_\bullet)$ is
{\em Cohen---Macaulay\/} if $\Gr^F_\bullet M$ is a Cohen---Macaulay module over
$\Gr^F_\bullet\cD$.
\end{defn}

A Cohen---Macaulay $\cD_X$-module $(M,F_\bullet)$ is holonomic iff
the dimension of $\Gr^F_\bullet M$ over $\Gr^F_\bullet\cD$ is equal
to $\dim X$.

\begin{lemma}
\label{lem:dual-to-restriction}
Suppose that $(M,F_\bullet)$ is a coherent holonomic filtered $\cD_X$-module
noncharacteristic with respect to $Y$.
Then $(M,F_\bullet)$ is holonomic Cohen---Macaulay at a point $y\in Y$
if and only if $(M,F_\bullet)_Y$ is.
\end{lemma}

\begin{pf}
We shall denote $(M,F_\bullet)_Y$ by $(M_Y,F)$.
Let $\dim X=n$, $\codim_X Y=d$.
Let $R=\Gr^F\cD_{X,y}$ and $R'=\Gr^F\cD_{Y,y}$, and
let $\m$ (respectively, $\m'$) denote the maximal ideal in $R$ (respectively,
$R'$) corresponding to the origin of $T^\ast_y X$ (respectively, $T^\ast_y Y$).
%generated by the image of the maximal ideal of $\cO_{Y,y}$ (respectively,
%$\cO_{X,y}$) and all elements of positive degree.

Both $R$ and $R'$ are graded rings, $\Gr^F M_y$ and $\Gr^F M_{Y,y}$ are
graded modules over them, and hence, the support of $\Gr^F M_y$ in $\Spec R$
corresponds to a homogeneous closed analytic subspace of $T^\ast U$ where
$U$ is a sufficiently small open neighborhood of $y$ in $X$, and
similarly for the support of $\Gr^F M_{Y,y}$  in $\Spec R'$.

We need to show that $\Gr^F M_y$ is Cohen---Macaulay of dimension $n$ over $R$
iff $\Gr^F M_{Y,y}$ is Cohen---Macaulay of dimension $n-d$ over $R'$.
The Cohen---Macaulay property of $\Gr^F M_y$
is equivalent to vanishing of $\Ext^j_R(\Gr^F M_y,R)$ for $j\ne n$;
 as the support of $\Ext^j_R(\Gr^F M_y,R)$
is homogeneous in $\Spec R$, this property holds at all points of
$\Spec R$ iff it holds at the origin of $T^\ast_y X$, i.e., at the
maximal ideal $\m$.
In other words, $\Gr^F M_y$ is Cohen---Macaulay of dimension $n$ over $R$ iff
$(\Gr^F M_y)_\m$ is Cohen---Macaulay of dimension $n$ over $R_\m$.
Similarly, $\Gr^F M_{Y,y}$ is a Cohen---Macaulay $R'$-module of
dimension $n-d$ iff
$(\Gr^F M_{Y,y})_{\m'}$ is a Cohen---Macaulay $R'_{\m'}$-module of
dimension $n-d$.

It follows that we need to show that $(\Gr^F M_y)_\m$
is Cohen---Macaulay of dimension $n$ over $R_\m$
iff $(\Gr^F M_{Y,y})_{\m'}$ is Cohen---Macaulay of dimension $n-d$
over $R'_{\m'}$.

Let $N=\Gr^F M_{Y,y}$.
Let $x_1,\dots,x_n$ be a local coordinate system in $X$ at $y$
such that $x_1,\dots,x_d$ is a system of local equations of $Y$ in $X$.
Since $Y$ is noncharacteristic, $x_1,\dots,x_d$ is a regular
$\Gr^F M_y$-sequence in $R$, and
$N\isomo\Gr^F M_y/(\sum_{l=1}^d x_l\Gr^F M_y)$.
Hence, $N$ is a module over $R/(\sum_{l=1}^d x_l R)$; its structure of an
$R'$-module comes from the embedding $R'\hookrightarrow R/(\sum_{l=1}^d x_l R)$.

Let $A$ and $A'$ be the quotients of $R/(\sum_{l=1}^d x_l R)$ and $R'$,
respectively, by the annihilators of $N$.  Then $A'\hookrightarrow A$;
since $\Gr^F M_Y$ is $\Gr^F\cD_Y$-coherent, $N$ is finite over $R'$ and
hence, $A$ is a finite $A'$-module.

Denote by $\tilde\m$ and $\tilde\m'$, respectively, the maximal ideals of
$A$ and $A'$ that correspond to the maximal ideals $\m$ and $\m'$ of $R$
and $R'$, respectively.
As $A$ is a finite $A'$-module, there are only finitely many ideals in $A$
lying over $\tilde\m'$, and clearly, $\tilde\m$ is one of them.
By a homogeneity argument, $\tilde\m$ is the only
ideal of $A$ lying over $\tilde\m'$.
Hence, $A_{\tilde\m}$ is a finite $A'_{\tilde\m'}$-module.

The localization $N_{\m'}=(\Gr^F M_{Y,y})_{\m'}$ of $N$ at $\m'$ as an
$R'$-module is the same as the localization $N_{\tilde\m'}$
of $N$ at $\tilde\m'$ as an
$A'$-module, and is isomorphic to the localization $N_{\tilde\m}$
of $N$ at $\tilde\m$ as an $A$-module (since $\tilde\m$ is the only
ideal of $A$ lying over $\tilde\m'$).

By~\cite{Se}, Ch.~IV, Proposition~12,
$\depth_{A'_{\tilde\m'}}N_{\tilde\m'}=\depth_{A_{\tilde\m}}N_{\tilde\m}$ and
$\dim_{A'_{\tilde\m'}}N_{\tilde\m'}=\dim_{A_{\tilde\m}}N_{\tilde\m}$.
Clearly, $\depth_{R'_{\m'}}N_{\m'}=\depth_{A'_{\tilde\m'}}N_{\tilde\m'}$
and $\dim_{R'_{\m'}}N_{\m'}=\dim_{A'_{\tilde\m'}}N_{\tilde\m'}$.

Since $x_1,\dots,x_d$ is a regular $\Gr^F M_y$-sequence in $R$, it is
a regular $(\Gr^F M_y)_\m$-sequence in $R_\m$.
We have $N_{\tilde\m}\isomo(\Gr^F M_y)_\m/(\sum_{l=1}^d x_l(\Gr^F M_y)_\m)$,
and hence,
$\depth_{A_{\tilde\m}}N_{\tilde\m}=\depth_{R_\m}(\Gr^F M_y)_\m-d$ and
$\dim_{A_{\tilde\m}}N_{\tilde\m}=\dim_{R_\m}(\Gr^F M_y)_\m-d$.

Altogether, we see that
$\depth_{R'_{\m'}}N_{\m'}=\depth_{R_\m}(\Gr^F M_y)_\m-d$ and
$\dim_{R'_{\m'}}N_{\m'}=\dim_{R_\m}(\Gr^F M_y)_\m-d$.
It follows that $\depth_{R_\m}(\Gr^F M_y)_\m=\dim_{R_\m}(\Gr^F M_y)_\m=n$
iff $\depth_{R'_{\m'}}N_{\m'}=\dim_{R'_{\m'}}N_{\m'}=n-d$, i.e.,
$(\Gr^F M_y)_\m$ is Cohen---Macaulay of dimension $n$ over $R_\m$ iff $N_{\m'}$
is Cohen---Macaulay of dimension $n-d$ over $R'_{\m'}$.
\end{pf}

\begin{remark}
\label{rem:right-exact}
It is not hard to see from the definitions that
the functors $\DR_X$ and $\aDR_X$ are right exact in the filtered sense:
if a filtered complex $(M^\bullet,F_\bullet)$ has the property that
$H^j\Gr^F_\bullet M^\bullet=0$ for $j>j_0$ then
$H^j\Gr^F_\bullet\DR_X(M^\bullet,F_\bullet)=0$ for $j>j_0$, and vice versa,
and the same holds for the functor of noncharacteristic restriction
$(\bullet)_Y$.
The functor $\aDR_X$ is also left exact in the similar sense.
\end{remark}

\section{Filtered perversity}
\label{sec:fil.perv}

In this section we introduce the notion of a {\em weakly filtered perverse}
differential complex; its meaning is that the complex is ``locally trivial''
in a certain filtered sense made precise below, and satisfies filtered
cohomology vanishing conditions which are similar to the cohomology vanishing
conditions for perverse sheaves.

The stratifications need to be defined only locally, which is made precise
by the notion of {\em stratified chart.}

This notion of {\em weak filtered perversity} is precisely
the assumption that we need to use; we call it {\em weak} because
we suspect that some stronger
property of ``local triviality'' along the strata will appear eventually.

We introduce also the notion of {\em coherent filtered perversity\/}
which is somewhat stronger than coherence together with weak
filtered perversity; we shall show in Proposition~\ref{cor:holon}
that it implies holonomicity of the corresponding $\cD_X$-module.

\subsection{Stratified charts}

\begin{defn}
A {\em stratified chart} $\cU$ on $X$ is the following
collection of data:
\begin{enumerate}
\item an open subset $U$ of $X$;
\item an analytic stratification of $U$;
\item for every point $x$ of any stratum $S$ of this stratification,
an open neighborhood $U_x$ of $x$ in $U$ and
an analytic submersion $\pi_x : U_x\to U_x\cap S$
which restricts
to the identity on $U_x\cap S$.
\end{enumerate}
\end{defn}

\subsection{The definition of filtered perversity}
In what follows we denote by $Y$ the fiber $\pi_x^{-1}(x)$.
Given a filtered differential complex $(\cM^\bullet,F_\bullet)$ on $X$, we use
the notation $\cF^j_p=H^j\Gr^F_p\aDR_X(\cM^\bullet,F_\bullet)$;
this is a sheaf of $\cO_X$-modules.

\begin{defn}
\label{cond:wfp}
A filtered differential complex $(\cM^\bullet,F_\bullet)$ on $X$ is called
{\em weakly filtered perverse} if $X$ can be covered by stratified charts
$\cU$ which satisfy
the following properties for every point $x$ of any stratum $S$ of $\cU$:
\begin{itemize}
\item[(i)]
for all $j$ and $p$, the sheaf $\cF^j_p$ has the property that
for all $i>0$, and for all
$y\in Y$, we have
$\Tor^{\cO_{X,y}}_i(\cF^j_{p,y},\cO_{Y,y}) = 0$;

\item[(ii)]
for all $j$, $p$, if
$\cF^j_{p,x}\otimes_{\cO_{X,x}}\cO_{Y,x} = 0$ then
$\cF^j_{p,x}=0$;

\item[(iii)] for all $p$, all $j<0$, we have
$H^j_{\{x\}}\Gr^F_p\DR_Y\left((\aDR_X(\cM^\bullet,F_\bullet))_Y\right) = 0$;

\item[(iv)] for all $p$, all $j>0$, we have $H^j\Gr^F_p\cM^\bullet = 0$.
\end{itemize}
We say that $(\cM^\bullet,F_\bullet)$ is {\em coherent filtered perverse\/}
if it is weakly filtered perverse, the complex
$\Gr^F_p\aDR_X(\cM^\bullet,F_\bullet)$ has bounded $\Gr^F_\bullet\cD_Y$-coherent
cohomology, and for any point $x$, $(\cM^\bullet,F_\bullet)$ is
noncharacteristic with respect to $Y$.
\end{defn}

Note that
property (i) means that the complex $\aDR_X(\cM^\bullet,F_\bullet)$
satisfies the condition \eqref{eqn:tor-nonchar}, and hence, the
noncharacteristic restriction $(\aDR_X(\cM^\bullet,F_\bullet))_Y$
which appears in (iii), is defined.

(We shall actually see that if $(\cM^\bullet,F_\bullet)$ is weakly filtered
perverse then $\aDR_X(\cM^\bullet,F_\bullet)$ is isomorphic to one filtered
$\cD_X$-module $(M,F_\bullet)$, and the condition of coherent filtered
perversity is equivalent to the condition that $(M,F_\bullet)$ coherent
holonomic Cohen-Macaulay.)

\subsection{Construction of stratified charts in the coherent case}
Suppose that $(\cM^\bullet,F_\bullet)$ is a filtered differential complex
such that $\Gr^F_\bullet\cM^\bullet$ has bounded coherent cohomology.
Then the condition (ii) is always satisfied.
We shall see here that under this assumption, there always exist
stratified charts satisfying also (i).

\begin{prop}
\label{prop:hom.shf.flat}
Suppose that $p:E\to X$ is a holomorphic vector bundle on $X$
and $\left\{\cF_i\right\}_I$ is a finite collection of homogeneous coherent
sheaves on $E$. Then at any point $x_0$ of
$X$ there exists a stratified chart $\cU$ such that for every
$x\in U$ and $i\in I$ the sheaf $\cF_i\vert_{p^{-1}U_x}$ is
$(\pi_xp)$-flat over $U_x\cap S$.
\end{prop}

\begin{pf}
Let $n=\dim X$.

We shall construct inductively a stratified chart $\cU^k$ containing
$x_0$ such that it satisfies the required flatness property at all
points $x$ of all the strata $S$ of codimension smaller than $k$.

We shall show how to construct $\cU^{k+1}$ once $\cU^k$ has been constructed.
Let $U^k$ be the open set containing $x_0$ which underlies $\cU^k$, and
let $X^k$ be the union of the closures of
all the strata of $\cU^k$ of codimension $k$.
Then $X^k$ is a closed analytic subset of $U^k$ of pure dimension $n-k$.
We choose an open polydisc $\Delta^n$ embedded in $U^k$ which
contains $x_0$ and such that there is a projection $q:\Delta^n\to\Delta^{n-k}$
with the property that the map
$q|_{X^k\cap \Delta^n}:X^k\cap \Delta^n\to\Delta^{n-k}$ is finite.

Consider the composite projection $qp:p^{-1}\Delta^n\to\Delta^{n-k}$.
Let $Z$ be the set of points in $p^{-1}\Delta^n\subset E$
where one of the sheaves $\cF_i$ is not $qp$-flat.
By Frisch's theorem on the openness of the flat locus
(\cite{F}, Theorem (IV,9) or \cite{BS}, Theorem V.4.5), $Z$ is a closed analytic
subset in $p^{-1}\Delta^n$ such that its image in $\Delta^{n-k}$ is negligible.
All the sheaves $\cF_i$ are homogeneous, and hence, $Z$ is homogeneous;
it follows that $p(Z)$ is a closed analytic subset of $\Delta^n$.
The intersection $p(Z)\cap X^k$ is negligible in $X^k\cap \Delta^n$
since its image
under $q$ is negligible and $q$ is finite on $X^k\cap \Delta^n$; it follows that
$p(Z)\cap X^k$ is a proper closed analytic subset of $X^k\cap \Delta^n$.

Construct $\cU^{k+1}$ as follows.
Take $U^{k+1}=\Delta^n$.
All the strata of $\cU^{k+1}$
of codimension less than $k$ are the intersections with
$U^{k+1}$ of the strata of $\cU^k$.
Any $(n-k)$-stratum $S'$ is obtained from a $(n-k)$-stratum $S$ of $\cU^k$
by intersecting with $U^{k+1}$ and then removing, first, all points where
$q|_{S\cap U^{k+1}}:S\cap U^{k+1}\to\Delta^{n-k}$ is ramified,
and second, the intersection with $p(Z)$.
The complement of these strata in $U^{k+1}$ has codimension at least
$k+1$; stratifying it, we complete the stratification of  $\cU^{k+1}$.

Clearly, the stata of codimension less than $k$ satisfy the required
flatness condition.
Let $S'$ be any $(n-k)$-stratum constructed as above.
At any point $x\in S'$ we take a neighborhood $U_x\subset U^{k+1}$
in such way that $q(U_x)=q(U_x\cap S')$ and
$q$ is an isomorphism on $U_x\cap S'$.
Take the projection $\pi_x:U_x\to U_x\cap S'$ such that $q\pi_x=q$, i.e.,
$\pi_x=(q|_{U_x\cap S'})^{-1}q$.
Then all the sheaves $\cF_i$ are $qp$-flat on $p^{-1}(U_x\cap S')$
since $S'$ does not intersect $p(Z)$, and hence, they are
$\pi_xp$-flat.
\end{pf}

\begin{cor}
\label{cor:strat}
Suppose that $(\cM_i^\bullet,F_\bullet)$ is a finite collection of
filtered differential complexes on $X$
such that $\Gr^F_\bullet\cM_i^\bullet$ have bounded coherent cohomology.
Then, locally
at any point of $X$ there exists a stratified chart such that
the properties (i) and (ii) of Definition~\ref{cond:wfp} hold for each
$(\cM_i^\bullet,F_\bullet)$.
\end{cor}

\begin{pf}
Consider the homogeneous coherent sheaves
$\left(H^j\Gr^F_\bullet\aDR(\cM_i^\bullet,F_\bullet)\right)^\sim$
on $T^*X$
obtained by localizing the corresponding $\Gr^F_\bullet\cD_X$-modules.

Proposition~\ref{prop:hom.shf.flat} yields a stratified chart satifying
the condition (i) for each
$(\cM_i^\bullet,F_\bullet)$; the condition (ii) is satisfied by coherence.
\end{pf}

\section{Filtered perverse complexes correspond to filtered
$\protect\cD$-modules}
\label{sec:fpc-to-dm}

\subsection{The main theorem}
Given a filtered complex $(M^\bullet,F_\bullet)$, the property that
$H^j\Gr^F_pM^\bullet =0$ for all $p$ and all $j\neq 0$, means that
$(M^\bullet,F_\bullet)$ is strict and $H^jM^\bullet = 0$ for $j\neq 0$.

Another formulation of the same property is that
$(M^\bullet,F_\bullet)$ is isomorphic to
$H^0(M^\bullet,F_\bullet)$ in the filtered derived category,
where $H^0(M^\bullet,F_\bullet)$ denotes $H^0M^\bullet$ equipped with the
induced filtration.

\begin{thm}
\label{thm:main}
Suppose that $(\cM^\bullet,F_\bullet)$ is weakly filtered perverse.
Then, for all $p$, all $j\neq 0$,
\begin{equation*}
H^j\Gr^F_p\aDR_X(\cM^\bullet,F_\bullet) = 0
\end{equation*}
Consequently the filtered complex
$\aDR_X(\cM^\bullet,F_\bullet)$ is
strict and isomorphic in the filtered derived category
to $H^0\aDR_X(\cM^\bullet,F_\bullet)$ equipped with the induced filtration.
\end{thm}

\begin{pf}
The statement is local so we may assume that $X=U$ in the definition of
weak filtered perversity, $x$ lies in the stratum $S$, $\pi_x: U_x\to S$ is
an analytic submersion which restricts to the identity on $S$, and
$Y=\pi_x^{-1}(x)$.

We are going to show by induction on $\codim S$
that the conclusion holds for the stalk of
$\Gr^F_p\aDR(\cM^\bullet,F_\bullet)$ at $x$. Thus we may assume that
the conclusion holds on the complement of the stratum $S$.

Condition (iv) of Definition \ref{cond:wfp} implies that, for all $p$,
\begin{equation}\label{van:above}
\text{$H^j\Gr^F_p\aDR_X(\cM^\bullet,F_\bullet) = 0$ for $j>0$}\ .
\end{equation}
By~\eqref{eqn:HGr(restr)} we have
%Condition (i) of Definition \ref{cond:wfp} implies that
\begin{equation}\label{iso:restr}
H^j\Gr^F_p\left((\aDR_X(\cM^\bullet,F_\bullet))_Y\right)\isomo
H^j\Gr^F_p\aDR_X(\cM^\bullet,F_\bullet)\otimes_{\cO_X}\omega_{Y/X}\ .
\end{equation}
%From \eqref{van:above} and \eqref{iso:restr} we obtain that
%\begin{equation}\label{van:above-on-Y}
%\text
%{$H^j\Gr^F_p\left((\aDR_X(\cM^\bullet,F_\bullet))_Y\right)= 0$ for $j>0$}\ .
%\end{equation}
The induction hypothesis and \eqref{iso:restr} imply that
\begin{equation}\label{van:on-Y-S}
\text{
$H^j\Gr^F_p\left((\aDR_X(\cM^\bullet,F_\bullet))_Y\right)
\vert_{Y\setminus\{x\}} = 0$
for $j\neq 0$}\ .
\end{equation}

Let $i:\{x\}\hookrightarrow Y$ be the embedding map.
Condition (iii) of Definition \ref{cond:wfp} implies that
$\R i^!\Gr^F_p\DR_Y\left((\aDR_X(\cM^\bullet,F_\bullet))_Y\right)$
is acyclic in negative degrees.
Here $\R i^!$ is the derived functor of the functor $i^!$ which assigns
to a sheaf $\cF$ its sections supported in $x$; if $\cF$ is an $\cO_Y$-module
or a $\cD_Y$-module then $\R i^!\cF$ is an an $\cO_{Y,x}$-module
or a $\cD_{Y,x}$-module, and $\R i^!$ commutes with the functors
$\Gr^F_p$ and $\DR_Y$.  Hence,
$\Gr^F_p\DR_Y\R i^!\left((\aDR_X(\cM^\bullet,F_\bullet))_Y\right)$
is acyclic in negative degrees; as $\aDR_Y$ is left exact
(Remark~\ref{rem:right-exact}), the complex
$\Gr^F_p\R i^!\left((\aDR_X(\cM^\bullet,F_\bullet))_Y\right)$
is acyclic in negative degrees, so that
\begin{equation}\label{van:loc-coh}
H^j_{\{x\}}\Gr^F_p\left((\aDR_X(\cM^\bullet,F_\bullet))_Y\right)= 0
\text{ for $j<0$ .}
\end{equation}
Examination of the long exact sequence in cohomology associated to the
inclusion $Y\setminus\{x\}\subset Y$ in the light of \eqref{van:on-Y-S} and
\eqref{van:loc-coh} shows that
\[
\text{$H^j\Gr^F_p\left((\aDR_X(\cM^\bullet,F_\bullet))_Y\right)= 0$
for $j<0$}\ .
\]
Together with \eqref{iso:restr} and the condition (ii) of Definition
\ref{cond:wfp} this shows that
\begin{equation}
\label{eqn:van:j<0}
\text{$H^j\Gr^F_p\aDR_X(\cM^\bullet,F_\bullet)_x = 0$ for $j<0$}\ .
\end{equation}
The statement of the Theorem is the combination of \eqref{van:above}
and \eqref{eqn:van:j<0}.
\end{pf}

\section{The coherent case}
\label{sec:coh.case}

In this section we study the property of filtered perversity of
a filtered differential complex $(\cM^\bullet,F_\bullet)$ under the
assumption of coherence of $H^\bullet\Gr^F_\bullet\cM^\bullet$; in particular,
this implies that the cohomology of $\Gr^F_\bullet\cM^\bullet$ is bounded.

By~\cite{S1}, (2.2.10.5), this is equivalent to the
$\Gr^F_\bullet\cD_X$-coherence of
$H^\bullet\Gr^F_\bullet\aDR(\cM^\bullet,F_\bullet)$; in case
$\aDR(\cM^\bullet,F_\bullet)$ is isomorphic to a single filtered
$\cD_X$-module, this property means that the module is $\cD_X$-coherent and
its filtration is good.

\subsection{Duality for coherent complexes}

The following technical lemma is a standard application of duality theory.

\begin{lemma}\label{lemma:duality}
Suppose that $X$ is a complex manifold of dimension $n$,
$L^\bullet$ is a bounded complex of coherent $\cO_X$-modules, and $x\in X$.
Then for each $j$, there is a nondegenerate pairing between
the spaces $H^{-j}_{\{x\}}L^\bullet$ and
$\left(H^j\R\shHom^\bullet_{\cO_X}(L^\bullet,\omega_X[n])\right)_x$;
the same is true for the spaces
$H^jL^\bullet_x$ and
$H^{-j}_{\{x\}}\R\shHom^\bullet_{\cO_X}(L^\bullet,\omega_X[n])$.
\end{lemma}

A nondegenerate pairing between two vector spaces
is a pairing that induces a monomorphism from each of them into
the (algebraic) dual of the other. (Actually, each of the vector spaces
can be given a topology so that they become topologically dual.
More precisely, the pairs of
spaces indicated in the Lemma, are strong dual to each other with respect to
certain natural FS and DFS topologies.
In case the complex $L^\bullet$ is zero except in one degree,
this statement
is a particular case of a theorem of Harvey: take $K=\{x\}$ in Theorem~5.12
of~\cite{ST}.
However, we do not need the topological duality;
all we need is that these vector spaces are either both zero or both nonzero.)

\begin{pf*}{Proof  of Lemma \protect\ref{lemma:duality}}
We shall establish the duality between the spaces
$H^{-j}_{\{x\}}L^\bullet$ and\break
$\left(H^j\R\shHom^\bullet_{\cO_X}(L^\bullet,\omega_X[n])\right)_x$;
the other duality would follow by substituting the dual complex
$\R\shHom^\bullet_{\cO_X}(L^\bullet,\omega_X[n])$ in place of $L^\bullet$.

Replacing $L^\bullet$ by its bounded free resolution in a neighborhood of $x$,
we may assume that all the sheaves $L^k$ are free.
Then $\R\shHom^\bullet_{\cO_X}(L^\bullet,\omega_X[n])$ is represented
by $\shHom^\bullet_{\cO_X}(L^\bullet,\omega_X[n])$, and
$\left(H^j\R\shHom^\bullet_{\cO_X}(L^\bullet,\omega_X[n])\right)_x\isomo
H^j\Hom^\bullet_{\cO_{X,x}}(L^\bullet_x,\omega_{X,x}[n])$.
The complex $\Hom^\bullet_{\cO_{X,x}}(L^\bullet_x,\omega_{X,x}[n])$
is a complex of free finitely generated $\cO_{X,x}$-modules.
Each of them has a canonical DFS topology and the differential is
continuous with respect to it; moreover, the image of the differential
is closed since it is closed with respect to the weaker topology of
coefficientwise convergence of formal power series (Theorem~6.3.5 of~\cite{H}).

Since each $L^k$ is free,
by a theorem of Martineau (\cite{ST}, Theorem 5.9) we have that $H^j_{\{x\}}L^k$
is zero for $j\ne n$, and $H^n_{\{x\}} L^k$ can be given a natural
Hausdorff FS topology in which it is a strong dual to
$\Hom_{\cO_{X,x}}(L_x^k,\omega_{X,x})$; in particular, it follows that
$H^{-j}_{\{x\}} L^\bullet\isomo H^{-j-n}(H^n_{\{x\}} L)^\bullet$ where
we denote
$(H^n_{\{x\}} L)^\bullet=
\{\dots\to H^n_{\{x\}} L^k\to H^n_{\{x\}} L^{k+1}\to\dots\}$.

The pairing between $H^n_{\{x\}} L^k$ and
$\Hom_{\cO_{X,x}}(L_x^k,\omega_{X,x})$, is given by the composition
of the multiplication
$H^n_{\{x\}}L^k\otimes\Hom_{\cO_{X,x}}(L_x^k,\omega_{X,x})\to
H^n_{\{x\}}\omega_X$
and the residue map $H^n_{\{x\}}\omega_X\to\bbC$, and hence, the complex
$(H^n_{\{x\}} L)^\bullet$ is the strong dual to the complex
$\Hom^\bullet_{\cO_{X,x}}(L_x^\bullet,\omega_{X,x}[n])$.
As the latter is a complex of DFS spaces with Hausdorff cohomology,
the former is a complex of FS spaces with Hausdorff cohomology
$H^{-j-n}(H^n_{\{x\}} L)^\bullet$ strong dual to
$\left(H^j\Hom^\bullet_{\cO_X}(L^\bullet,\omega_X[n])\right)_x$.
This yields a nondegenerate pairing between $H^{-j}_{\{x\}} L^\bullet$
and $\left(H^j\Hom^\bullet_{\cO_X}(L^\bullet,\omega_X[n])\right)_x$.
\end{pf*}

\subsection{Holonomicity}
\begin{prop}
\label{cor:holon}
Suppose that $(\cM^\bullet,F_\bullet)$ is a coherent filtered perverse
complex on a complex manifold $X$. Then
$\aDR(\cM^\bullet,F_\bullet)$ is isomorphic in the filtered derived category
to a filtered holonomic Cohen-Macauley $\cD_X$-module.
\end{prop}

\begin{pf}
The question is local and we need to prove it in a neighborhood of
any point $x\in X$.

The point $x$ is covered by a stratified chart $\cU$ satisfying properties
(i)--(iv) of Definition~\ref{cond:wfp}; we keep the notation introduced there.

Our assumptions imply that $\aDR_X(\cM^\bullet,F_\bullet)$ is isomorphic to a
coherent filtered $\cD_X$-module; we shall denote this module by
$(M,F_\bullet)$.

We argue by induction by the codimension of the stratum $S$ containing
$x$; by the inductive assumption, we may assume that $(M,F_\bullet)$
is holonomic Cohen-Macauley in the complement to $S$.

By Lemma~\ref{lem:dual-to-restriction},
this implies that $(M,F_\bullet)_Y$ is holonomic
Cohen-Macauley everywhere on $Y$ except possibly at $x$; hence, it
is holonomic.

Property (iii) implies that $H^j_{\{x\}}\Gr^F_p\DR_Y((M,F_\bullet)_Y)=0$
if $j<0$.
By Lemma~\ref{lemma:duality} this yields
$\left(H^j\R\shHom^\bullet_{\cO_Y}
(\Gr^F_p\DR_Y((M,F_\bullet)_Y),\omega_Y[\dim Y])\right)_x=0$
if $j>0$.
By \S\ref{subsec:duality}
this implies that $H^j\Gr^F_p\DR_Y\bbD((M,F_\bullet)_Y)=0$ if $j>0$,
and by right exactness of $\aDR_Y$ we get
$H^j\Gr^F_p\bbD((M,F_\bullet)_Y)=0$ if $j>0$.

Since $(M,F_\bullet)_Y$ is a filtered $\cD_Y$-module, we have
$H^j\Gr^F_p\DR_Y((M,F_\bullet)_Y)=0$ for $j>0$.
By Lemma~\ref{lemma:duality} and \S\ref{subsec:duality}, this implies
$H^j_{\{x\}}\Gr^F_p\DR_Y\bbD((M,F_\bullet)_Y)=0$ for $j<0$.
So we get $H^j_{\{x\}}\Gr^F_p\bbD((M,F_\bullet)_Y)=0$ for $j<0$
by the left exactness of $\aDR_Y$.
The long exact sequence of the inclusion $Y\setminus\{x\}\subset Y$
(cf.\ the proof of Theorem~\ref{thm:main}) yields
$H^j\Gr^F_p\bbD((M,F_\bullet)_Y)=0$ if $j<0$.
(Actually, this vanishing also follows from the vanishing --- see, for
example,~\cite{Borel}, V.2.2.2 --- of
$\Ext^i_{\Gr^F\cD_{Y,x}}(\Gr^F M_{Y,x},\Gr^F\cD_{Y,x})$ for $i<d$ where
$d$ is the codimension of the support of $\Gr^F M_{Y,x}$ in
$\Spec\Gr^F\cD_{Y,x}$; in our case $d=\dim Y$.)

Altogether, we see that $H^j\Gr^F_p\bbD((M,F_\bullet)_Y)=0$ if $j\ne 0$,
i.e., the filtered complex $\bbD((M,F_\bullet)_Y)$ is isomorphic to one
filtered module.  Hence, $(M,F_\bullet)_Y$ is Cohen-Macaulay at $x$.

It follows by Lemma~\ref{lem:dual-to-restriction}
that $(M,F_\bullet)$ is also Cohen-Macaulay at $x$.
\end{pf}

\subsection{The converse to the Main Theorem}

\begin{thm}
\label{thm:converse}
If a coherent filtered $\cD_X$-module $(M,F_\bullet)$ is holonomic
and Cohen-Macauley, then $\DR_X(M,F_\bullet)$ is coherent filtered perverse.
\end{thm}

\begin{pf}
Consider a point $x\in X$.

The coherence of $(M,F_\bullet)$ implies that
$H^\bullet\Gr^F_\bullet\DR(M,F_\bullet)$ is $\cO_X$-coherent.
Consequently, Corollary~\ref{cor:strat}
implies that in a neighborhood of $x$ there
exists a stratified chart $\cU$ such that $\DR(M,F_\bullet)$
satisfies properties (i) and (ii) with respect to it.
We shall keep the notation of Definition~\ref{cond:wfp}.

As $M$ is holonomic, we may assume that that the characterisic
variety $\Ch(M)$ is contained in the union of the conormal bundles
to the strata of the Whitney stratification that underlies $\cU$.
It follows that the projection $(Y\times_X T^\ast X)\cap\Ch(M)\to T^\ast Y$
is finite (it is even an embedding), and hence, $(M,F_\bullet)_Y$
is a coherent filtered $\cD_Y$-module.
This implies that $Y$ is noncharacteristic with respect to $(M,F_\bullet)$.

The property (iv) of Definition~\ref{cond:wfp}
at $x$ is satisfied by $\DR(M,F_\bullet)$
since it is satisfied by the de Rham complex of any filtered $\cD_X$-module.

Let us show the property (iii) at $x$ with respect to this stratified chart.

By Lemma~\ref{lem:dual-to-restriction},
$(M,F_\bullet)_Y$ is holonomic Cohen-Macaulay at $x$.
Hence, the complex
$\bbD((M,F_\bullet)_Y)$ is isomorphic to a filtered $\cD_Y$-module,
and consequently, the complex $\DR_Y\bbD((M,F_\bullet)_Y)$ satisfies (iv)
at $x$:
\[
H^j\Gr^F_\bullet\DR_Y\bbD((M,F_\bullet)_Y)=0\text{ \ at $x$ for all $j>0$.}
\]
By \S\ref{subsec:duality} and Lemma~\ref{lemma:duality} we get
\[
H^j_{\{x\}}\Gr^F_\bullet\DR_Y((M,F_\bullet)_Y)=0\text{ \ for all $j<0$.}
\]
This is the property (iii) at $x$ for the filtered complex $\DR_X(M,F_\bullet)$.
\end{pf}

\section{An application to {\protect\lt} cohomology}
\label{sec:appl}

\renewcommand{\Fl}{{F_\bullet}}

In this section we give an application to our results and show that
in the situation of \cite{KK2} (and under the assumption of filtered
perversity of the \lt\ complex) there is a local filtered isomorphism
(in the sense of derived category) between the $L^2$ complex
and the de Rham complex of
the $\cD$-module that underlies the corresponding pure Hodge module.

\subsection{}
\label{subsec:l2}
Let $X$ denote a K\"ahler manifold of dimension $n$, let
$j:{\xcirc}\hookrightarrow X$ be the inclusion map of the complement of
a divisor with normal crossings, and
$\bbE = (\bbE_{\bbQ},(\cO_{\xcirc}\otimes_{\bbQ}\bbE_{\bbQ},\Fl))$
a quasiunipotent  polarised variation of pure Hodge structure of weight
$w$ on $\xcirc$.

Let $(N,F_\bullet)$ denote the filtered $\cD_X$-module underlying the
polarizable Hodge module~\cite{S1} which restricts to
$(\omega_{\xcirc}\otimes_{\bbQ}\bbE_{\bbQ},\Fl)$ on $\xcirc$.

Let $(\ltcx,\Fl)$ denote
the {\lt}-complex with coefficients in $\bbE$ constructed using the Hodge
inner product in the fibers of $\bbE$ and a certain complete metric $\eta$ on
$\xcirc$ as in~\cite{KK2}, \cite{CKS}; to keep up with our degree conventions,
we shall assume that the grading of $(\ltcx,\Fl)$
is chosen in such way that ${\cal A}^i_{(2)}({\Bbb E})$ contains forms of
degree $i+n$.

As the metric $\eta$ satisfies $\eta>C\eta_X$ locally in a neighborhood
of any point of $X$, where $\eta_X$ is
the metric on $X$ and $C$ a suitable positive constant,
the holomorphic forms on $X$ are bounded in the pointwise
norm with respect to $\eta$ and the {\lt}-complex is an $\cO_X$-module:
if $\omega$ is a section of $\ltcx$ and $f$ is a holomorphic function
then $f\omega$ is also a section of $\ltcx$ since both $f\omega$
and $d(f\omega)=df\wedge\omega+fd\omega$ are \lt.

We shall assume here that $(\ltcx,\Fl)$ is weakly filtered
perverse, and so is every direct summand of it
(in the sense of derived category).
By Theorem~\ref{thm:main} this implies that the complex
$\aDR(\ltcx,\Fl)$ is strict and isomorphic in the filtered derived category
to its zeroeth cohomology with the induced filtration.

By~\cite{KK2} and \cite{CKS}, $\ltcx$ is isomorphic in the derived category of
complexes
of sheaves on $X$ to the intersection complex with coefficients in $\bbE$.
We shall assume, moreover, that for any cross-section $Y$
appearing in the definition of weak filtered perversity, the
complex
$\DR_Y\left(\left(\aDR_X({\cal A}^{\bullet}_{(2)}(X,{\Bbb E}))\right)_Y\right)$
is isomorphic
to the intersection cohomology complex on $Y$ with the coefficients in
${\Bbb E}|_{Y\cap\xcirc}$
in the derived category of complexes (without filtration).
(One would even expect that the filtered complex
$\DR_Y\left(\left(\aDR_X({\cal A}^{\bullet}_{(2)}(X,{\Bbb E}),F_\bullet)
\right)_Y\right)$
is isomorphic in the filtered derived category to
$({\cal A}^{\bullet}_{(2)}(Y,{\Bbb E}|_{Y\cap\xcirc}),F_\bullet)$.)

In case $X$ is compact, both complexes of global sections
$\Gamma(X,\ltcx)$ and $\Gamma(X,\DR(N))$ are strict (\cite{KK2}, \cite{S1}),
and their cohomology have pure Hodge structures.
Their cohomology groups are isomorphic (\cite{KK2}) together with the Hodge
filtrations (\cite{S2}, p.~294).  Here we strengthen these results and
show the isomorphism at the level of sheaves (without the assumption of
compactness), in the filtered derived categories:

\begin{prop}
\label{thm:applic}
Assume that 
\begin{enumerate}
\item any direct summand of
$({\cal A}^{\bullet}_{(2)}(X,{\Bbb E}),F_\bullet)$
in the filtered derived category of filtered differential complexes
(in particular,
$({\cal A}^{\bullet}_{(2)}(X,{\Bbb E}),F_\bullet)$ itself)
is weakly filtered perverse;
\item for any cross-section $Y$
as in Definition \ref{cond:wfp} (weak filtered perversity), the
complex
$\DR_Y\left(\left(\aDR_X({\cal A}^{\bullet}_{(2)}(X,{\Bbb E}))\right)_Y\right)$
is isomorphic
to the intersection cohomology complex on $Y$ with the coefficients in
${\Bbb E}|_{Y\cap\xcirc}$.
\end{enumerate}
Then the filtered differential complexes $\DR_X(M,F_\bullet)$ and
$({\cal A}^{\bullet}_{(2)}(X,{\Bbb E}),F_\bullet)$ are isomorphic
in the filtered derived category.
Equivalently, the filtered $\cD_X$-modules
$(M,F_\bullet)$ and $H^0\aDR(\ltcx,F_\bullet)$ are isomorphic.
\end{prop}

\begin{pf}
By Remark 3.15 of~\cite{S2} (the idea actually going back
to~\cite{KK}), there is a direct sum decompostion in the
derived category of filtered differential complexes
\begin{equation*}
(\ltcx,\Fl)\isomo \DR_X(N,F_\bullet)\oplus (\cM^\bullet,F_\bullet)
\end{equation*}
and we need to show that the second summand is trivial.

Our assumptions imply that
\begin{enumerate}
\item
$(\cM^\bullet,F_\bullet)$ is weakly filtered
perverse, therefore, by Theorem~\ref{thm:main}, isomorphic to
$\DR_X(M,F_\bullet)$ where $(M,F_\bullet)$ is a filtered $\cD_X$-module;
\item
for any cross-section $Y$ as in Definition \ref{cond:wfp}, the complex
$\DR_Y((\aDR_X\cM^\bullet)_Y)\isomo\DR_Y(M_Y)$ is acyclic.
\end{enumerate}

The weak filtered perversity of $(\cM^\bullet,F_\bullet)$
implies that $X$ is covered by stratified charts satisfying properties
(i)--(iv) of Definition \ref{cond:wfp}. It is sufficient to show that the
intersection of the support of $H^\bullet\Gr^F_\bullet(\cM^\bullet,F_\bullet)$
with any of the charts is empty.

Assume to the contrary and consider a stratified chart and a point
$x\in\supp H^\bullet\Gr^F_\bullet(\cM^\bullet,F_\bullet)$ which lies
on a stratum which is maximal among those which have a nonempty
intersection with $\supp H^\bullet\Gr^F_\bullet(\cM^\bullet,F_\bullet)$.

Let $Y$ denote the cross-section at $x$. Then $\supp M_Y\subseteq\{x\}$.
In addition we have
$H^j\Gr^F_\bullet\DR_Y((M,F_\bullet)_Y)=0$ for $j>0$ by the right exactness
of $\DR_Y$, and for $j<0$ by property (iii) observing that
$H^j\Gr^F_\bullet\DR_Y((M,F_\bullet)_Y)\isomo
H^j_{\{x\}}\Gr^F_\bullet\DR_Y((M,F_\bullet)_Y)$
since $\Supp\Gr^F_\bullet\DR_Y((M,F_\bullet)_Y)\subseteq\{x\}$.

It follows that the filtered complex $\DR_Y((M,F_\bullet)_Y)$ is strict.
Since the complex $\DR_Y(M_Y)$ is acyclic it follows that
the complex $\DR_Y((M,F_\bullet)_Y)$ is filtered acyclic, and hence,
$M_Y$ is trivial.
The property (ii) of the weak filtered perversity shows that $M$
is trivial at $x$ (cf.\ the proof of Theorem~\ref{thm:main}),
which contradicts our assumption. Hence, $M$ is trivial.
\end{pf}

\begin{cor}
In the assumptions of Proposition~\ref{thm:applic} the sheaves
$H^\bullet\Gr^F_\bullet\ltcx$ (the \lt-$\dbar$-cohomology, see
below) are coherent.
\end{cor}

\subsection{Remarks on the \lt-$\overline\partial$-cohomology of
a singular variety}
\label{rem:bound-cond}
Let $(\cM^\bullet,F_\bullet)$ be the \lt-complex of a singular subvariety
$Z$ (the complex of sheaves of forms with locally summable coefficients
on the nonsingular part $Z^\circ$ of $Z$ which
are \lt\ together with their differentials near all points of $Z$,
both smooth and singular).

This is the sheafification of the presheaf assigning to each open set
$U\subset Z$ the domain of the maximal closed extension of the
differential $d$ on the Hilbert space of the \lt-forms on $U\cap Z^\circ$;
there is another flavor of the \lt-complex constructed in a similar
way by sheafification of the {\em minimal\/} closed extension, see details
in~\cite{Y}, \S2.3.

The complex $\Gr^F_{-p}\cM^\bullet$ consists of $(p,q)$-forms
$\omega^{pq}$ on $Z^\circ$ (for any $q$)
which have the following properties:
$\omega^{pq}$ and $\dbar\omega^{pq}$ are \lt\ and moreover,
there exists a form $\omega=\omega^{pq}+\omega^{p+1,q-1}+\dots$ such
that both $\omega$ and $d\omega$ are \lt.

The
differential in the complex $\Gr^F_p\cM^\bullet$ is the operator $\dbar$.
Let us {\em assume\/} that this is a closed extension of $\dbar$.
(To be precise, this means that the sections of this sheaf over an
open set $U\subset Z$ form a closed extension of $\dbar$ in the
Fr\'echet space of forms on $U\cap Z^\circ$ which are \lt\ locally
in a neighborhood of any point of $U$; the topology on this Fr\'echet
space is given by the seminorms $\|\bullet\|_K$ where $K$ is a relatively
compact open subset
of $U$; for a form $\omega$, the value $\|\omega\|_K$ is the \lt\ norm
of $\omega$ on $K$.)

In such case the complex $\Gr^F_p\cM^\bullet$ can be viewed as an
``ideal boundary condition'' (the notion due to J. Cheeger)
for the operator $\dbar$ at the singularities of $Z$;
this complex contains the minimal closed extension of $\dbar$
and is contained in the maximal one (their sheafifications
can be defined in a way similar to those of the operator $d$).

If the operator $d$ on \lt\ forms has the property that its minimal
extension coincides with the maximal one (this is called the {\em $L^2$
Stokes property\/}~\cite{C}, \cite{Y}; it is known for conical
singularities~\cite{C} and seems from~\cite{O3}
to be a reasonable conjecture in general) then, under our assumptions,
it is not hard to see that the boundary condition
of $\Gr^F_{-p}\cM^\bullet$ is dual to the boundary condition
of $\Gr^F_{-p'}\cM^\bullet$ if $p+p'=\dim Z$. (More precisely, this
means the following.
For any open $U\subset Z$, the dual to
the Fr\'echet space of forms which are locally
L2 on $U$ --- with the topology described above --- can be identified
by the pairing
$<\omega,\phi>=\int_{U\cap Z^\circ}\omega\wedge\phi$ with the DF space
of forms $\phi$ on $U\cap Z^\circ$ such that the closure of $\supp\phi$
in $U$ is compact.
The $L^2$ Stokes property means that the adjoint of maximal extension of
$d$ in the first space is, up to sign,
the maximal extension of $d$ in the second one.
The differential $\dbar$ on the sections of $\Gr^F_{-p}\cM^\bullet$ on $U$
is an unbounded operator on the subspace of $(p,q)$-forms of the above
Fr\'echet space; similarly, the differential $\dbar$ on the sections of
$\Gr^F_{-p'}\cM^\bullet$ with compact support is an unbounded operator
on the subspace of $(p',q)$-forms of the above DF space.
The duality between the boundary conditions means that these two operators
are adjoint up to sign.)

It is easy to see that under our assumptions,
the differential in $\Gr^F_{-\dim Z}\cM^\bullet$ is the maximal closed
extension of $\dbar$; it follows by duality that the differential in
$\Gr^F_0\cM^\bullet$ is the minimal closed extension of $\dbar$.

In case the metric on $Z^\circ$ is complete (e.g., Saper metric),
the minimal closed extension of $\dbar$ is known to coincide with the
maximal one, and hence, $\Gr^F_{-p}\cM^\bullet$ is just the domain
of $\dbar$ with any of the boundary conditions.

In case the metric is incomplete (e.g., the restriction of the Fubini---Studi
metric on the projective space to $Z^\circ$), it is known that the
minimal closed extension of $\dbar$ may be different from the
maximal one~\cite{P}.  The results of~\cite{PS} and~\cite{FH}
suggest that in case
$p=0$ the ``correct'' boundary condition for $\dbar$ is the minimal
(Dirichlet) one, and in case
$p=\dim Z$ the ``correct'' boundary condition for $\dbar$ is the maximal
(Neumann) one.

This suggests that under our assumptions,
the complex $\Gr^F_{-p}\cM^\bullet$ is the most natural
boundary condition for the operator $\dbar$, and its cohomology
$H^j\Gr^F_{-p}\cM^\bullet$ is the most natural notion of the
\lt-$\dbar$-cohomology sheaves.

\end{document}